\documentclass{amsart}

\usepackage[english]{babel}
\usepackage{ae,aecompl}
\usepackage{amsmath,amsfonts,amssymb,amsthm,latexsym,stmaryrd} 
\usepackage{amscd,graphicx}
\usepackage[dvipsnames]{color}
\usepackage[all]{xy}

\newcommand{\nc}[2]{\newcommand{#1}{#2}}
\newcommand{\rnc}[2]{\renewcommand{#1}{#2}}
\rnc{\theequation}{\thesection.\arabic{equation}}
\def\note#1{{}}

\nc{\beq}{\begin{equation}}
\nc{\eeq}{\end{equation}}
\rnc{\[}{\beq}
\rnc{\]}{\eeq}
\nc{\ba}{\begin{array}}
\nc{\ea}{\end{array}}
\nc{\bea}{\begin{eqnarray}}
\nc{\beas}{\begin{eqnarray*}}
\nc{\eeas}{\end{eqnarray*}}
\nc{\eea}{\end{eqnarray}}
\nc{\be}{\begin{enumerate}}
\nc{\ee}{\end{enumerate}}
\nc{\bd}{\begin{diagram}}
\nc{\ed}{\end{diagram}}
\nc{\bi}{\begin{itemize}}
\nc{\ei}{\end{itemize}}
\nc{\bpr}{\begin{proposition}}
\nc{\bth}{\begin{theorem}}
\nc{\ble}{\begin{lemma}}
\nc{\bco}{\begin{corollary}}
\nc{\bre}{\begin{remark}\em}
\nc{\bex}{\begin{example}\em}
\nc{\bde}{\begin{definition}}
\nc{\ede}{\end{definition}}
\nc{\epr}{\end{proposition}}
\nc{\ethe}{\end{theorem}}
\nc{\ele}{\end{lemma}}
\nc{\eco}{\end{corollary}}
\nc{\ere}{\hfill\mbox{$\Diamond$}\end{remark} }
\nc{\eex}{\hfill\mbox{$\Diamond$}\end{example}}
\nc{\bpf}{{\it Proof.~~}}
\nc{\epf}{\hfill\mbox{$\square$}\vspace*{3mm}}
\nc{\hsp}{\hspace*}
\nc{\vsp}{\vspace*}
\nc{\lN}{{\ell}^2(\bN)}
\nc{\lZ}{{\ell}^2(\bZ)}
\newcommand{\wegdamit}[1]{}

\nc{\ot}{\otimes}
\nc{\te}{\!\ot\!}
\nc{\bmlp}{\mbox{\boldmath$\left(\right.$}}
\nc{\bmrp}{\mbox{\boldmath$\left.\right)$}}
\nc{\LAblp}{\mbox{\LARGE\boldmath$($}}
\nc{\LAbrp}{\mbox{\LARGE\boldmath$)$}}
\nc{\Lblp}{\mbox{\Large\boldmath$($}}
\nc{\Lbrp}{\mbox{\Large\boldmath$)$}}
\nc{\lblp}{\mbox{\large\boldmath$($}}
\nc{\lbrp}{\mbox{\large\boldmath$)$}}
\nc{\blp}{\mbox{\boldmath$($}}
\nc{\brp}{\mbox{\boldmath$)$}}
\nc{\LAlp}{\mbox{\LARGE $($}}
\nc{\LArp}{\mbox{\LARGE $)$}}
\nc{\Llp}{\mbox{\Large $($}}
\nc{\Lrp}{\mbox{\Large $)$}}
\nc{\llp}{\mbox{\large $($}}
\nc{\lrp}{\mbox{\large $)$}}
\nc{\lbc}{\mbox{\Large\boldmath$,$}}
\nc{\lc}{\mbox{\Large$,$}}
\nc{\Lall}{\mbox{\Large$\forall\;$}}
\nc{\bc}{\mbox{\boldmath$,$}}
\nc{\ra}{\rightarrow}
\nc{\ci}{\circ}
\nc{\cc}{\!\ci\!}
\nc{\lra}{\longrightarrow}
\nc{\imp}{\Rightarrow}
\rnc{\iff}{\Leftrightarrow}
\nc{\inc}{\mbox{$\,\subseteq\;$}}
\rnc{\subset}{\inc}

\nc{\0}{\sb{(0)}}
\nc{\3}{\sb{(3)}}
\nc{\4}{\sb{(4)}}
\nc{\5}{\sb{(5)}}
\nc{\6}{\sb{(6)}}
\nc{\7}{\sb{(7)}}

\def\Tr{{\rm Tr}}

\def\<{\langle}
\def\>{\rangle}

\def\id{\mbox{$\mathop{\mbox{\rm id}}$}}

\nc{\spp}{\mbox{${\mathcal S}{\mathcal P}(P)$}}
\nc{\ob}{\mbox{$\Omega\sp{1}\! (\! B)$}}
\nc{\op}{\mbox{$\Omega\sp{1}\! (\! P)$}}
\nc{\oa}{\mbox{$\Omega\sp{1}\! (\! A)$}}
\nc{\dr}{\mbox{$\Delta_{R}$}}
\nc{\dsr}{\mbox{$\Delta_{\Omega^1P}$}}
\nc{\ad}{\mbox{$\mathop{\mbox{\rm Ad}}_R$}}
\nc{\as}{\mbox{$A(S^3\sb s)$}}
\nc{\bs}{\mbox{$A(S^2\sb s)$}}
\nc{\slc}{\mbox{$A(SL(2,\C))$}}
\nc{\suq}{\mbox{$\cO(SU_q(2))$}}
\nc{\CS}{\mbox{$C(S^1)$}}
\nc{\tc}{\widetilde{can}}
\def\slq{\mbox{$\cO(SL_q(2))$}}

\nc{\p}{{\rm pr}}

\nc{\barot}{\bar \otimes}

\nc{\ha}{\mbox{$\alpha$}}
\nc{\hb}{\mbox{$\beta$}}
\nc{\hg}{\mbox{$\gamma$}}
\nc{\hd}{\mbox{$\delta$}}
\nc{\he}{\mbox{$\varepsilon$}}
\nc{\hz}{\mbox{$\zeta$}}
\nc{\hs}{\mbox{$\sigma$}}
\nc{\hk}{\mbox{$\kappa$}}
\nc{\hm}{\mbox{$\mu$}}
\nc{\hn}{\mbox{$\nu$}}
\nc{\hl}{\mbox{$\lambda$}}  
\nc{\hG}{\mbox{$\Gamma$}}
\nc{\hD}{\mbox{$\Delta$}}
\nc{\hT}{\mbox{$\Theta$}}
\nc{\ho}{\mbox{$\omega$}}
\nc{\hO}{\mbox{$\Omega$}}
\nc{\hp}{\mbox{$\pi$}}
\nc{\hP}{\mbox{$\Pi$}}

\nc{\qpb}{quantum principal bundle}

\def\H{{\Bbb H}}
\def\C{{\Bbb C}}

\def\R{{\Bbb R}}

\def\cO{{\mathcal O}}

\def\cD{{\mathcal D}}

\def\cK{{\mathcal K}}

\def\pr{{\rm pr}}

\def\tc{\tilde{c}}

\def\a{\alpha} \def\b{\beta}  
\def\p{\psi}

\def\p{\psi}

\def\bC{\mathbb{C}}
\def\bN{\mathbb{N}} 

\def\bR{\mathbb{R}}  
 
\def\bZ{\mathbb{Z}}

\def\sC{\mathcal{C}} 
\def\sO{\mathcal{O}} \def\sT{\mathcal{T}}

\def\ex{\mathrm{e}}

 \def\id{\mathrm{id}}

\def\op{\mathrm{op}}

\def\Tr{\mathrm{Tr}}

\def\be{\begin{equation}}
\def\ee{\end{equation}}

\def\ot{\otimes}

\def\hx{\hat{x}}

\def\hy{\hat{y}}
\def\hatyap{{\hy'}_\alpha}
\def\hz{\hat{z}}

\def\a{\alpha}
\def\b{\beta}
\def\as{\alpha^\ast}
\def\bs{\beta^\ast}
\def\slq2{SL_q (2)}

\def\half{{\frac{1}{2}}}

\def\CSUq2{C( SU_q (2))}
\def\ASUq2{A( SU_q (2))}
\def\OSUq2{\sO ( SU_q (2))}

\def\oa{\overline{\alpha}}

\def\ob{\overline{\beta}}

\newtheorem{theorem}{Theorem}[section]
\newtheorem{proposition}[theorem]{Proposition}
\newtheorem{lemma}[theorem]{Lemma}
\newtheorem{corollary}[theorem]{Corollary}

\theoremstyle{definition}
\newtheorem{definition}[theorem]{Definition}
\newtheorem{example}[theorem]{Example}
\newtheorem{remark}[theorem]{Remark}

\begin{document}

\title{Index pairings for pullbacks of C*-algebras}
\vspace*{-15mm}

\author{Ludwik D\k abrowski}
\address{SISSA (Scuola Internazionale Superiore di Studi Avanzati)\\
Via Bonomea 265, 34136 Trieste, Italy\\
E-mail: dabrow@sissa.it}

\author{Tom Hadfield}
\address{
G-Research\\ 
Whittington House, 19-30 Alfred Place, London WC1E 7EA, United Kingdom\\
E-mail: Thomas.Daniel.Hadfield@gmail.com}

\author{Piotr M.~Hajac}
\address{Instytut Matematyczny, Polska Akademia Nauk\\
ul.\ \'Sniadeckich 8, 00-956 Warszawa, Poland\\
Katedra Metod Matematycznych Fizyki, Uniwersytet Warszawski\\
ul.\ Ho\.za 74, 00-682 Warszawa, Poland\\
E-mail: pmh@impan.pl}

\author{Rainer Matthes}
\address{Fachgruppe Mathematik,
Fakult\"{a}t f\"{u}r Physikalische Technik und Informatik,\\
Wests\"{a}chsische Hochschule Zwickau\\
Dr.-Friedrichs-Ring 2A, 08056 Zwickau, Germany\\
E-mail:  rainer.matthes@fh-zwickau.de}

\author{Elmar Wagner}
\address{Instituto de F\'isica y Matem\'aticas, Universidad Michoacana\\
Edificio C-3, Cd. Universitaria, C.\,P.\,58040 Morelia, Mich., Mexico\\
E-mail: elmar@ifm.umich.mx}

\maketitle
\vspace*{-5mm}

\mbox{ }\hfill {\it Dedicated to S.\ L.\ Woronowicz on the occasion of 
his 70th birthday}
\vspace*{3mm} 

\begin{abstract}
In this overview, we study how to reduce the index pairing for a fibre-product 
C*-algebra to the index pairing for the C*-algebra over which the fibre 
product is taken.
As an example we analyze the case of suspensions and apply 
it to noncommutative instanton bundles of arbitrary charges over 
the suspension of quantum deformations of the 3-sphere.
\end{abstract}

\section{Introduction} 

The aim of this paper is to review and put in a general context a number 
of results concerning the index computations for different
types of noncommutative 
instanton vector bundles. By the Serre-Swan theorem,
the category of finite-dimensional vector bundles 
over compact Hausdorff spaces is equivalent to the category
of finitely generated projective modules over algebras
of all continuous functions on the corresponding
compact Hausdorff spaces. The latter category makes perfect sense
for arbitrary algebras, so that finitely generated projective
modules over noncommutative algebras are considered as noncommutative 
vector bundles.

 In this paper, we prove some general theorems concerning
finitely generated projective modules over one-surjective
pullbacks of 
C*-algebras.
Then we specialize to non-reduced suspensions of C*-algebras viewed
as the pullbacks of cones of C*-algebras. Choosing them
to be noncommutative deformations of the algebra of all continuous 
functions on the 3-sphere, we can construct idempotents of
noncommutative instanton vector bundles for an arbitrary 
``winding number" (cf.~\cite{atiyah})
and determine their index
pairing with an appropriate K-homology class (``instanton charges").

Our key idea is to take advantage of the compatibility of the index 
pairing with  the connecting homomorphisms of 
the Mayer-Vietoris six-term exact sequences for K-theory and K-homology,
which is a manifestation of
 the associativity  of the Kasparov product.
This allows us to reduce the computation of an index pairing  
from non-reduced suspensions of C*-algebras to the C*-algebras
themselves.
Geometrically speaking, we determine instanton charges by shifting
calculations from a quantum 4-sphere to its equator quantum 3-sphere,
 which considerably simplifies the task.

Herein we focus on idempotents of noncommutative instanton vector
bundles coming from a particular 
complex-parameter deformation of the algebra
$C(S^3)$ of all continuous functions on the unital 3-dimensional
sphere (cf.~\cite{BCT1,lpr}). 
The real and unitary cases of this deformation
correspond respectively to Woronowicz's C*-algebra $C(SU_q(2))$
\cite{wo} of the quantum  group $SU_q(2)$ 
and Matsumoto's C*-algebra $C(S^3_\theta)$ \cite{matsumoto}
of the quantum sphere $S^3_\theta$ built from noncommutative 
tori~\cite{r-ma90}.

\section{The Mayer-Vietoris six-term exact sequence}
\label{Mayer-Vietoris}

\noindent
Let $B_1 \stackrel{\pi_1}{\longrightarrow} A 
\stackrel{\pi_2}{\longleftarrow} B_2$
 be homomorphisms
of C*-algebras.
We denote by $B$ the C*-algebra fiber product (pullback
C*-algebra) \cite{Pedersen}
\[\label{pull}
B_1 \underset{A}{\times}  B_2
:=\{ (b_1,b_2)\in B_1 \times B_2 ~|~ \pi_1(b_1)= \pi_2(b_2)\}.
\]
With the obvious projection maps $\pr_1 : B \to B_1$ 
and $\pr_2 : B \to B_2$,
there is
the commutative diagram
\begin{equation}\label{B_is_fiber_product} 
\mbox{$\xymatrix@=5mm{& & B \ar[lld]_{\pr_1}
\ar[rrd]^{\pr_2}& &\\
B_1 \ar[drr]_{\pi_1}& & & &B_2\,, \ar[dll]^{\pi_2}\\
&& A&&}$}
\end{equation}
and $B$ enjoys a universal property with respect to all such
commutative diagrams. 

To use Mayer-Vietoris arguments in K-theory, we need
to assume that one of the defining homomorphisms is surjective.
In this paper we choose  $\pi_2$ 
to be surjective. To go beyond K-theory and use  results of
\cite[Section~1.2.3]{bm}, we also assume that $\pi_2$ is semi-split.
This is always the case for surjections onto nuclear C*-algebras,
which covers all our examples.
First, there exist the following Mayer-Vietoris 
six-term exact sequences both in K-theory and K-homology: 
\begin{equation}\label{mv}
     	\begin{CD}
	{K_0 (B)} @ >{({\pr_1}_\ast , {\pr_2}_\ast)}>> 
	{K_0 ( B_1 ) \oplus K_0 (B_2)} @ >{{\pi_2}_\ast - {\pi_1}_\ast} >> {K_0 ({A} )} @ . @  .\\ 
 	@ A{\partial_{10}}AA @ . @ VV
     {\partial_{01}}
V @ .\\
	{K_1 ({A} ) } @ <{{\pi_2}_\ast - {\pi_1}_\ast} << 
	{K_1 (B_1 ) \oplus K_1 (B_2)} @ <{({\pr_1}_\ast , {\pr_2}_\ast)} << {K_1 (B),} @ . @ .\\
     	\end{CD}
\end{equation}
\
 \begin{equation}\label{mv_khom}
     	\begin{CD}
	{K^0 (B)} @ <{ \pr_1^\ast  + \pr_2^\ast}<< 
	{K^0 ( B_1 ) \oplus K^0 (B_2)} @ <{( - \pi_1^\ast , \pi_2^\ast )} << {K^0 ({A} )} @ 
	. @  .\\ 
 	@ V{\delta_{01}}VV @ . @ AA    {\delta_{10}} A @ .\\
	{K^1 ({A} ) } @ >{( - \pi_1^\ast, \pi_2^\ast )} >> 
	{K^1 (B_1 ) \oplus K^1 (B_2)} @ >{\pr_1^\ast + \pr_2^\ast } >> {K^1 (B)} @ . @ . .\\
     	\end{CD}
\end{equation}
Moreover, our paper  applies the fundamental result 
that the Mayer-Vietoris connecting homomorphisms  
are compatible
with the index pairing \cite[Section~1.2.3]{bm}. 
These connecting homomorphisms are determined
uniquely up to a sign, so that their compatibility with the index pairing
is also determined up to a sign. Herein we use an explicit
construction of $\partial_{10}$ and choose $\delta_{01}$ to
be such that
for any ${\bf x}  \in K_1({A})$ and 
${\bf w}\in K^0 (B)$
 \begin{equation}\label{compatibility}
\boxed{< \partial_{10}  ({\bf x} ),   {\bf w}  > \,=\, 
< {\bf x} , \delta_{01} ( {\bf w} ) >} \, .
\end{equation}

\subsection{Odd-to-even connecting homomorphism}
\label{subsection_Bass}

Given a ring homomorphism 
 \mbox{$\stackrel{f}{R \to S}$} and a left $R$-module $E$, 
define $f_\ast E$ to be the left $S$-module $S \otimes_R E$. 
 Then there is a canonical $R$-linear map $f_{\ast} : E \to f_\ast E$ 
given by $f_\ast(e) := 1 \otimes_R  e$. 
   It is clear that, if $E$ is free over $R$ with basis 
$\{ \hx_\a \}_\alpha$, then $f_\ast E$ is free over $S$ with basis 
$\{ x_\a := f_\ast ( \hx_\a) = 1 \otimes_R \hx_\a \}_\alpha$.
Next, following \cite{milnor}, we take 
 modules $E_1$ and $E_2$ over the rings $B_1$ and $B_2$
respectively, and observe that the Abelian group $E_1 \times E_2$ is a 
left module over the pullback ring \eqref{pull} via 
$b \cdot ( \xi_1 , \xi_2 ) := ( \pr_1 (b) \xi_1 , \pr_2 (b) \xi_2 )$. 
Given  an isomorphism $h : {\pi_1}_\ast E_1 \cong {\pi_2}_\ast E_2$, we 
can now define 
\[\label{bmod}
 E := M( E_1, E_2 ,h):=
\{( \xi_1 , \xi_2)\in E_1\times E_2\;|\;
(h \circ {\pi_1}_\ast) ( \xi_1 ) = {\pi_2}_\ast ( \xi_2)\}.
\] 
 
Milnor's theorem \cite[Theorem~2.1]{milnor} asserts that,
if $E_1$ and $E_2$ are finitely generated projective over $B_1$ and
$B_2$ respectively, then $E$ is finitely generated projective
over~$B$. In particular, 
choosing $E_1$ and $E_2$ to be the free modules $B_1^N$
and $B^N_2$ respectively, \mbox{$N\in\bN\setminus\{0\}$},
 this beautiful construction 
allows Milnor to
define a connecting homomorphism \cite[p.\ 28]{milnor}
for the Mayer-Vietoris long
exact sequence in algebraic K-theory. Herein we prove
the aforementioned special case of \cite[Theorem~2.1]{milnor}
supplementing it with an explicit formula for an idempotent
matrix representing the finitely generated projective module~$E$.
Our explicit formula agrees with the formula provided by Ranicki
in his unpublished notes~\cite{Ran}.\footnote{We are grateful
to Ulrich Kr\"ahmer for making us aware of this reference.}

\begin{theorem}\label{milnor_proj}
Let $B$ be the pullback ring \eqref{pull} and $E$ be the left
$B$-module~\eqref{bmod}. Assume that $E_i\cong B_i^N$, $i\in\{1,2\}$,
$N\in\bN\setminus\{0\}$, as left modules, and take 
$a:=(a_{\a\b})\in GL_N(A)$
to be the matrix implementing an isomorphism 
$h\colon{\pi_1}_\ast(E_1)\to{\pi_2}_\ast(E_2)$. Then,
if $\pi_2$ is surjective,
$E \cong B^{2N} p$ as a left module, where 
 \begin{equation}\label{P}
 p:= 
  \left[
  \begin{array}{cc}
  (1, c(2-dc)d) & (0, c(2-dc)(1-dc)) \cr
  (0, (1-dc)d) & (0, (1-dc)^2 )
  \end{array}
  \right]
 \end{equation}
is an idempotent matrix in $M_{2N}(B)$ such that 
all entries of the sub-matrices $c,d\in M_N(B_2)$ satisfy
$\pi_2(c_{\a\b})=a_{\a\b}$  
 and $\pi_2(d_{\a\b})=(a^{-1})_{\a\b}$ for all $\alpha$ and $\beta$.
 \end{theorem}
 \begin{proof} 
Let  $\{ \hx_\a \}_{1 \leq \a \leq N}$ and  
$\{ \hy_\b \}_{1 \leq \b \leq N}$ be bases of $E_1 \cong B_1^N$ 
and $E_2 \cong B_2^N$ respectively.
Similarly, let  $\{ x_\a := {\pi_1}_\ast (\hx_\a )\}_\alpha$ and 
$\{ y_\b := {\pi_2}_\ast ( \hy_\b ) \}_\beta$ be the corresponding bases
for ${\pi_1}_\ast E_1$ and ${\pi_2}_\ast E_2$ respectively.
 Both of these modules are isomorphic 
to $A^N$. Therefore,
the isomorphism $h : {\pi_1}_\ast E_1 \to {\pi_2}_\ast E_2$ is given 
by $h( x_\a ) =: \sum_\beta a_{\a \b} y_\b$ for some unique
$(a_{\a \b }) \in GL_N ( A)$. Now, the matrix $c$ might be
invertible or non-invertible. If it is, we quickly conclude
that $E\cong B^N$. Without this assumption we need a doubling 
construction to prove~\eqref{P}. 

\emph{Special case.} Assume that $c$ is invertible. 
   For each $\a$, define $\hatyap := \sum_\b c_{\a \b} \hy_\b$. 
   Since $c$ is invertible, $\{ \hatyap \}_\a$ is a basis of $B_2^N$. 
  Also, as $(h \circ {\pi_1}_\ast) ( \hx_\a ) = h( x_\a ) 
= \sum_\b a_{\a \b} y_\b = {\pi_2}_\ast ( \hatyap)$,
    it follows that $\hz_\a := ( \hx_\a , \hatyap ) \in E:= M( E_1 , E_2 , h)$ 
for all $\a$. Finally,
   it is  straightforward to check that 
$\{ \hz_\a \}_{1 \leq \a \leq N}$ 
span $E$ and are linearly independent, so that $E \cong B^N$ 
via the basis $\{ \hz_\a \}_\a$.
   
   \emph{General case.} Let us use the matrix $a^{-1}$ to
define the isomorphism 
\[
g \colon {\pi_1}_\ast E_1 \to {\pi_2}_\ast E_2,\quad
g ( x_\a ) 
:= \sum_\b (a^{-1})_{\a \b} y_\b,\quad \a\in\{1,\cdots, N\}.
\]
 Then $h \oplus g$ defines an isomorphism 
 ${\pi_1}_\ast ( E_1 \oplus E_1) \to {\pi_2}_\ast 
( E_2 \oplus E_2)$. Next,
as in Whitehead's lemma, we write
 \[
\left[
  \begin{array}{cc}
  a & 0 \cr
  0 & a^{-1}
  \end{array}
  \right]
 = 
 \left[
  \begin{array}{cc}
  1 & a \cr
  0 & 1
  \end{array}
  \right]
  \left[
  \begin{array}{cc}
  1 & 0 \cr
  - a^{-1} & 1
  \end{array}
  \right]
 \left[
  \begin{array}{cc}
  1 & a \cr
  0 & 1
  \end{array}
  \right]
 \left[
  \begin{array}{cc}
  0 & -1 \cr
  1 & 0
  \end{array}
  \right]
  \in M_{2N} (A),
   \]
where $1$ means the $N\times N$-identity matrix~$I_N$.
    This lifts to
   \begin{align}\label{liftc}
   \sC &:= 
    \left[
  \begin{array}{cc}
  \sC^{11} & \sC^{12} \cr
  \sC^{21} & \sC^{22}
  \end{array}
  \right]\nonumber\\
 &:= 
 \left[
  \begin{array}{cc}
  1 & c \cr
  0 & 1
  \end{array}
  \right]
  \left[
  \begin{array}{cc}
  1 & 0 \cr
  - d & 1
  \end{array}
  \right]
 \left[
  \begin{array}{cc}
  1 & c \cr
  0 & 1
  \end{array}
  \right]
 \left[
  \begin{array}{cc}
  0 & -1 \cr
  1 & 0
  \end{array}
  \right]\nonumber\\
  &\phantom{:}= 
  \left[
  \begin{array}{cc}
  (2-cd)c & cd-1 \cr
  1-dc & d
  \end{array}
  \right]
    \in M_{2N} (B_2),
    \end{align}
 whose    inverse is
 \[\label{liftd}
 \cD := 
  \left[
  \begin{array}{cc}
  \cD^{11} & \cD^{12} \cr
  \cD^{21} & \cD^{22}
  \end{array}
  \right]
 := 
 \left[
  \begin{array}{cc}
  d & 1-dc \cr
  cd-1 & (2-cd)c
  \end{array}
  \right]
  \in M_{2N} (B_2).
 \]
This allows us to proceed as in the special case. First, however,
we need to introduce more notation. We put 
\[
\hat{x}_\a^L:=(\hat{x}_\a,0),\quad \hat{x}_\a^R:=(0,\hat{x}_\a),
\quad\a\in\{1,\cdots,N\},
\]
 to obtain 
$\{ \hx_\a^L, \hx_\b^R \}_{ 1 \leq \a, \b \leq N}$ as a basis
of $E_1\oplus E_1$, and 
\[
\hat{y}_\a^L:=(\hat{y}_\a,0),\quad 
\hat{y}_\a^R:=(0,\hat{y}_\a),\quad
\a\in\{1,\cdots,N\},
\]
 to obtain 
$\{ \hy_\a^L, \hy_\b^R \}_{ 1 \leq \a, \b \leq N}$ as a basis
of $E_2\oplus E_2$.
 Furthermore, we put $(\sC^{ij}_{\a \b} ) :=\sC^{ij}\in M_N ( B_2)$
for all $i,j\in\{1,2\}$.
  Finally, we    
  define
  \[\label{solve}
  \hat{u}_\a := \sum_{\mu}\Big( \sC^{11}_{\a \mu} \; {\hy}_{\mu}^L 
+ \sC^{12}_{\a \mu} \; {\hy}_{\mu}^R \Big)\;, \quad
  \hat{v}_\b := \sum_{\mu} \Big( \sC^{21}_{\b \mu} \; {\hy}_{\mu}^L 
+  \sC^{22}_{\b \mu} \; {\hy}_{\mu}^R\Big)\,.
\]
Now we can argue as in the special case to conclude that
 $M(E_1 \oplus E_1 , E_2 \oplus E_2 , h \oplus g)
\cong B^{2N}$ through the basis 
$\{ ( \hx_\a^L, \hat{u}_\a ), ( \hx_\b^R , \hat{v}_\b ) \}_{1 \leq \a, \b \leq N}$. 
On the other hand, we have the isomorphism
\begin{align}
M( E_1 , E_2 , h) \oplus M( E_1 , E_2 , g) 
\,&\longrightarrow\, 
M( E_1 \oplus E_1 , E_2 \oplus E_2 , h \oplus g ),
\nonumber\\
( (e_1, e_2) , (f_1 , f_2 )) 
&\longmapsto ( ( e_1 , f_1 ) , ( e_2 , f_2 )).
\end{align}
Consequently, $E := M( E_1 , E_2 , h)$ is a direct summand of a 
finitely generated free module, so that it is a finitely generated
 projective module, as claimed.

We proceed now to computing the idempotent matrix~$p \in M_{2N} (B)$.
Composing the projection and inclusion
\[
M( E_1 \oplus E_1 , E_2 \oplus E_2 , h \oplus g) 
\stackrel{\pi}{\ra} M( E_1 , E_2 , h) \stackrel{i}{\hookrightarrow}
M( E_1 \oplus E_1 , E_2 \oplus E_2 , h \oplus g)
\]
yields an idempotent: $(i \circ \pi)^2 = i \circ \pi$. 
To determine the matrix $p$ of
$i \circ \pi$ with respect to the basis 
$\{ ( \hx_\a^L, \hat{u}_\a ), ( \hx_\b^R , \hat{v}_\b ) \}_{1 \leq \a, \b \leq N}$,
we solve~\eqref{solve} for ${\hy}_{\mu}^L$ with the help of~\eqref{liftd}
\[
{\hy}_{\mu}^L= \sum_{\nu}\Big(    
\cD^{11}_{\mu \nu}\hat{u}_{\nu}+
\cD^{12}_{\mu \nu}\hat{v}_{\nu}\Big),
\]
and compute
\begin{align}
(i \circ \pi) ( \hx_\a^L , \hat{u}_\a ) 
&= (i \circ \pi) \Big( \hx_\a^L\;\text{\large,}\; \sum_{\mu}\big(
 \sC^{11}_{\a \mu} \; {\hy}_{\mu}^L 
+  \sC^{12}_{\a \mu} \; {\hy}_{\mu}^R\big) \Big) \\
&= \Big( \hx_\a^L \;\text{\large,}\;  \sum_{\mu} \sC^{11}_{\a \mu} \; {\hy}_{\mu}^L\Big)\nonumber\\
&=\Big( \hx_\a^L \;\text{\large,}\;  
\sum_{\mu,\nu} \sC^{11}_{\a \mu}\big(    
\cD^{11}_{\mu \nu}\hat{u}_{\nu}+
\cD^{12}_{\mu \nu}\hat{v}_{\nu}
\big)\Big)\nonumber\\
&= \sum_{\nu}\Big((I_N)_{\a\nu}\;\text{\large,}\;
(\sC^{11}\cD^{11})_{\a \nu}\Big)( \hx_{\nu}^L , \hat{u}_{\nu} )
+\sum_{\nu}\Big(0\;\text{\large,}\;
(\sC^{11}\cD^{12})_{\a \nu}\Big)( \hx_{\nu}^R , \hat{v}_{\nu} ),
\nonumber
\end{align}
\begin{align}
(i \circ \pi) (\hx_\b^R , \hat{v}_\b ) 
&= (i \circ \pi) \Big(\hx_\b^R \;\text{\large,}\; 
 \sum_{\mu} \big( \sC^{21}_{\b \mu} \; {\hy}_{\mu}^L 
+  \sC^{22}_{\b \mu} \; {\hy}_{\mu}^R\big) \Big)\\
&= \Big(0\;\text{\large,}\; \sum_{\mu} \; 
\sC^{21}_{\b \mu} \; {\hy}_{\mu}^L \Big)\nonumber\\
&= \Big(0\;\text{\large,}\; \sum_{\mu,\nu} \; 
\sC^{21}_{\b \mu} \big(    
\cD^{11}_{\mu \nu}\hat{u}_{\nu}+
\cD^{12}_{\mu \nu}\hat{v}_{\nu}
\big) \Big)\nonumber\\
&= \sum_{\nu}\Big(0\;\text{\large,}\;
(\sC^{21}\cD^{11})_{\b \nu}\Big)( \hx_{\nu}^L , \hat{u}_{\nu} )
+\sum_{\nu}\Big(0\;\text{\large,}\;
(\sC^{21}\cD^{12})_{\b \nu}\Big)( \hx_{\nu}^R , \hat{v}_{\nu} ).
\nonumber
\end{align}
Combining this with \eqref{liftc} and \eqref{liftd} yields
\[
p = 
\left[\!
  \begin{array}{cc}
  (1,\sC^{11}\cD^{11} ) & (0, \sC^{11}\cD^{12}) \cr
  (0, \sC^{21}\cD^{11}) & (0, \sC^{21}\cD^{12} )
  \end{array}
 \! \right]
=
 \left[\!
  \begin{array}{cc}
  (1, c(2-dc)d) & (0, c(2-dc)(1-dc)) \cr
  (0, (1-dc)d) & (0, (1-dc)^2 )
  \end{array}
 \! \right],
  \]
as desired.
\end{proof}

Furthermore, it is useful to observe
  that the idempotent $p$ given by \eqref{P} 
can be written as 
\begin{equation}\label{array}
\boxed{
 p = X Y, \quad
 X:=\left[
\begin{array}{c}
(1, c(2-dc)) \cr
(0, 1-dc)
\end{array}
\right], \quad
Y:= \left[ (1,d), (0, 1-dc)\right]
}\,.
\end{equation}
Since $Y X =(1, 1)\in B$, we  immediately see that
$p^2=p$.
 Note that the components of $X$ and $Y$ are in 
$B_1\times B_2$ but not necessarily  in $B$, which allows the non-freeness
of the projective module given by $p$. 

Finally, let us use the above theorem to obtain an explicit formula
for a connecting homomorphism. Herein we adapt Milnor's construction
to unital C*-algebras as proposed by Nigel Higson 
in~\cite[Section~0.4]{hrz}.
\begin{theorem}\label{milnorhigson} 
Let $B$ be the pullback C*-algebra
\eqref{pull} of unital C*-ho\-mo\-mor\-phisms\linebreak
\mbox{$B_1 \stackrel{\pi_1}{\longrightarrow} A 
\stackrel{\pi_2}{\longleftarrow} B_2$}. 
Assume that $\pi_2$ is surjective.
Also, let $p\in M_{2N}(B)$ be an idempotent assigned to 
$a\in GL_N(A)$ by~\eqref{P}, and let $I_N\in M_N(B)$
denote the $N\times N$-identity matrix.
Then the formula 
\begin{equation}\label{partial}
\partial_{10} \colon K_1(A)\ni [a] \longmapsto [p] - [I_N ]\in K_0(B) 
 \end{equation}
defines a connecting homomorphism of the Mayer-Vietoris six-term exact 
sequence~\eqref{mv}.
\end{theorem}

\subsection{Even-to-odd connecting homomorphism}
\label{subsection_Bott}

Using the Bott periodicity to factorize $\partial_{01}$ into
\[
K_0(A)\stackrel{\cong}{\longrightarrow} 
K_1(C_0(\mathbb{R})\otimes A)\longrightarrow K_1(B),
\]
one can directly and 
explicitly compute  the even-to-odd connecting homomorphism
of the Mayer-Vietoris six-term exact sequence~\cite{bhr}.
Alternatively, one can simply combine \cite[Theorem~1.18]{bm}
with \cite[Section~9.3.2]{blackadar} to obtain:
\begin{theorem}
Let $p\in M_N(A)$ be a projection, and  
$Q$ be a self-adjoint lifting of $p$ to~$B_2$, i.e.\
 $\pi_2(Q)=p$ and $Q^*=Q$. Also, let $I_N\in M_N(B)$
denote the $N\times N$-identity matrix.
Then  the formula  
\begin{equation}\label{baumhajac}
\boxed{
\partial_{01}\colon K_0(A)\ni([p])\longmapsto
[(I_N,\ex^{2\pi i Q})]\in K_1(B)
}
\end{equation}
defines a connecting homomorphism of the Mayer-Vietoris six-term exact 
sequence~\eqref{mv}.
\end{theorem}

\subsection{Suspension case}\label{section_suspension}

\noindent
In this section, we assume that $B$ is the fibre product of two cones of $A$
over $A$: 
\begin{equation}
B_1 := \{ f \in  C( [0,\mbox{$\frac{1}{2}$}] , {A}) \; | \; f(0) \in \bC \},  
\qquad
B_2 := \{ f \in  C( [\mbox{$\half$},1] , {A}) \; | \;  f(1) \in \bC  \},
\end{equation}
with $\pi_1 : B_1 \to A$ and $\pi_2 : B_2 \to A$  both given by
 $f \mapsto f(\half)$.
Clearly,
$B$ can be identified with the non-reduced suspension of $A$
(as in algebraic topology) via the isomorphism:
\begin{equation}\label{via}
\Sigma A := 
\{f \in  C( [0,1] , {A}) \; | \; f(0), f(1) \in \bC \}
\ni f \longmapsto (f|_{[0,\half ]}, f|_{[\half ,1]})\in B.
\end{equation}
In this particular setting, the idempotent \eqref{P}
can be given more explicitly.
Given  $a \in GL_{N} (A)$, we lift it to 
$c :=\psi a \in M_N (B_2)$,
where 
$\psi \in C([\mbox{$\frac{1}{2}$}, 1])$
satisfies $\psi(\mbox{$\frac{1}{2}$}) = 1$ and $\psi(1) = 0$. The inverse $a^{-1}$ we lift to
 $d:= \psi a^{-1}\in M_N (B_2)$.
(In what follows we shall further restrict the choice of $\psi$
to fit some known examples.) Plugging in $c$ and $d$ to \eqref{array} yields 
\begin{align} \begin{split}
\label{xysusp}
X&=\left[
\begin{array}{c}
(1, \psi  (2-\psi^2)a)\cr
(0,(1-\psi^2))
\end{array}
\right]\in M_{2N\times N}(B_1\times B_2),    \\
Y&=\big[(1,\psi a^{-1})\,,\,(0,(1-\psi^2) )\big]
\in M_{N\times 2N}(B_1\times B_2). 
\end{split}
\end{align}

Next, since cones are contractible,
\begin{equation}
K_j (B_i ) = 
\left\{
\begin{array}{cc}
\bZ & \text{for} \;\; j=0,\cr
0 & \text{for} \;\; j=1,
\end{array}
\right.\quad i=1,2 \, .
\end{equation}
Consequently, the six-term exact sequence
 (\ref{mv}) reduces to the exact sequence
\begin{equation} \label{5lem}
\xymatrixcolsep{4pc}\xymatrix{
0 \to K_1 (A) \stackrel{\partial_{10}\ }{\longrightarrow} K_0 (B) 
\ar[r]^{\phantom{.}\hspace{30pt}({\pr_1}_\ast , {\pr_2}_\ast )} 
& \bZ \oplus \bZ 
\ar[r]^{\phantom{.}\hspace{-40pt}{\pi_2}_\ast - \,{\pi_1}_\ast } 
& K_0 ( A) \stackrel{\partial_{01}\ }{\longrightarrow}
K_1 (B) \to 0\, .
}
\end{equation}
Assume furthermore that $K_0 (A) \cong\bZ$ via $[1]\mapsto 1$, 
and that $K_1(A) \cong \bZ$ via $[\mathsf{u}]\mapsto 1$, for some 
$\mathsf{u}\in U_N(A)$.
It follows that the map 
$
{\pi_2}_\ast - {\pi_1}_\ast:\bZ \oplus \bZ\rightarrow  \bZ
$
is surjective, so that $K_1 (B) = 0$.
On the other hand, as the kernel of ${\pi_2}_\ast 
- {\pi_1}_\ast \cong \bZ$,
by standard homological algebra we conclude that
$K_0 (B)\cong \bZ^2$ with one generator given by $1\in B$ and
the other given by $\partial_{10}([\mathsf{u}])$ defined in~\eqref{partial}.

Now, by Rosenberg and Schochet's Universal Coefficient 
Theorem~\cite{rs87}, 
if the K-groups $K_j $ are free, then $K^j  \cong K_j $ 
(as Abelian groups). 
Thus we obtain
\begin{equation}
K^j (A)= 
\left\{
\begin{array}{cc}
\!\!\!\bZ &\!\!\!\text{for }  j=0, \cr
\!\!\!\bZ &\!\!\!\text{for }  j=1,
\end{array}
\right.\!\!\!\quad
K^j (B_i)= 
\left\{
\begin{array}{cc}
\!\!\!\bZ &\!\!\!\text{for }  j=0, \cr
\!\!\! 0 &\!\!\!\text{for }  j=1,
\end{array}
\right.\!\!\!\quad
K^j (B)= 
\left\{
\begin{array}{cc}
\!\!\!\bZ^2 &\!\!\!\text{for }  j=0, \cr
\!\!\! 0 &\!\!\!\text{for }  j=1.
\end{array}
\right.\!\!\!
\end{equation}
Combining it with the K-homology six-term exact 
sequence~\eqref{mv_khom}, we infer that the connecting
homomorphism
$\delta_{01} : K^0 (B) \to K^1 ( A)$ is surjective.
Hence, for any ${\bf z}\in K^1 (A)$ we can choose $ {\bf w}\in K^0 (B)$
such that  $\delta_{01}({\bf w})={\bf z}$, and use the compatibility
\eqref{compatibility}
of the index pairing with the Mayer-Vietoris six-term exact sequences
to compute
\begin{equation}\label{main1}
< \partial_{10} ( [\mathsf{u}^n]) ,  
{\bf w}  >\, = \,< [\mathsf{u}^n] ,  {\bf z}  > 
\,=\, n< [\mathsf{u}] , {\bf z} >,\quad n\in\bZ\, .
\end{equation}
Remembering that the index pairing with a free module is always zero,
we may now conclude that computing a non-trivial index pairing for an 
idempotent $p_n$ corresponding to the unitary $\mathsf{u}^n$, 
$n\in\bZ$, via
\eqref{P} boils down to finding an odd Fredholm module whose pairing
with $\mathsf{u}$ equals one:
\begin{equation}\label{main2}
< [p_n],  {\bf w}  > 
\,=\, n< [\mathsf{u}] , {\bf z} > .
\end{equation}

Finally, let us remark that for the usual suspension of $A$, given by
\begin{equation}
SA := 
\{ \; f \in  C( [0,1] , {A}) \; | \; f(0)= f(1)=0 \},
\end{equation} 
there is also a method to construct an idempotent 
$\widetilde{p}\in M_{2N} (SA)$
in terms of a unitary 
$\mathsf{u}\in U_{N} (A)$ \cite[p.~139]{wegge-olsen}.
This method applied to our non-reduced suspension $\Sigma A$ of $A$
 yields $\widetilde{p}: = T T^*$, where 
\begin{equation}\label{woproj}
T:=
\left[
\begin{array}{c}
\cos^2(\mbox{$\frac{\pi t}{2}$})+\sin^2(\mbox{$\frac{\pi t}{2}$})
\, \mathsf{u}  \cr 
\cos(\mbox{$\frac{\pi t}{2}$})\sin(\mbox{$\frac{\pi t}{2}$})
(\mathsf{u}^{-1} -1)
\end{array}
\right],\quad\text{and $T^*$ is its Hermitian conjugate.}
\end{equation}

\section{Noncommutative instanton vector bundles of an arbitrary charge}
\label{section_nc_inst}

In this section, we focus on the case when $A$ is the C*-algebra
of a quantum or classical~$S^3$, so that the foregoing
assumptions about the K-groups of $A$ are satisfied. 
To begin with, we take $A:=C(S^3)$. Then $B\cong\Sigma C(S^3)\cong
C(S^4)$. We know that $K_1(C(S^3))$
is generated by the fundamental representation of $SU(2)$ identified with $S^3$.
Explicitly, $u\in U_2(C(S^3))$ can be described as the continuous map
\begin{eqnarray}\begin{split}
&S^3\ni {x} = ( x_0, x_1, x_2, x_3 ) \,\stackrel{u}{\longmapsto}\,  
\left[ 
\begin{array}{cc}
 a,  & -\bar b\cr
b, &  \bar a 
\end{array}
\right]\in SU_2(\bC)\, , &\\ 
\label{pauli}
&\text{where} \quad 
a:=x_0+ i x_3,\quad b:=x_1 + i x_2\,. &
\end{split}
 \end{eqnarray}
Likewise, we can introduce the coordinate functions
\begin{equation}  \label{uclassic}
S^3\ni {x} \,\stackrel{\alpha}{\longmapsto} a\in\bC,
\qquad
S^3\ni {x}  \,\stackrel{\beta}{\longmapsto} b\in\bC,
\qquad
u:=\left[ 
\begin{array}{cc}
 \a,  & -\b^*\cr
\b, &  \a^* 
\end{array}
\right].
\end{equation}
In what follows, we will also view $u$ as a function on 
the space of unit quaternions.

As before, we can take $p_n$ to be an idempotent corresponding
via \eqref{P} to the unitary $u^n$, $n\in\bZ$. Then, by the Serre-Swan
Theorem, the finitely generated projective module $(C(S^4))^4 \, p_n$
can be identified with the module of all continuous sections of 
the uniquely determined (up to an isomorphism) rank 2 complex
vector bundle
$E_n\to S^4$. These vector bundles are commonly referred to as 
\emph{instanton} vector bundles.
One can prove that the ``winding number" $n\in\bZ$ defining the
vector bundles $E_n$ coincides with
a $K_0$-invariant, which we call the charge of an instanton.

Using homotopy and coordinate changes, we  obtain many
different  but $K_0$-equivalent explicit descriptions of the idempotents
$p_n$, $n\in\bZ$. First, recall that identifying
 $\R^4$ with the space of quaternions $\H$ 
we can define the  unit 4-sphere $S^4\subset\bR^5$ as
\begin{equation}
S^4:=\{(\tau, h)\in\bR\times\H\;|\;\tau^2 + \| h\|^2=1\}\,.
\end{equation}
 Another way is to present $S^4$ is via the non-reduced suspension of $S^3$.
Since the latter can be identified with the space of quaternions of norm~1,
we obtain\\[-30pt]
\begin{align}
&\label{cong}\\
&S^4\cong\left\{(t, g)\in [0,1]\times\H\;{\big|}\; \| g\|=1 \right\}/
\!\sim\,,\quad (0,g)\sim(0,g'),
\quad (1,g)\sim(1,g'),\;\forall\; g,g'\in\H\, .\nonumber
\end{align}
Finally, by removing the north or the south pole from $S^4$, we can
cover it by two quaternionic charts via stereographic projections.
Denoting one of these stereographic coordinates
 by $z\in\H$, we obtain the following
relations between the above defined coordinate systems:
\begin{equation}
\label{transf}
g=\frac{h}{\sqrt{1-\tau^2}}, \  \,t= \frac{1+\tau}{2},\qquad
z=\sqrt{\frac{t}{1-t}}\, g, \qquad
z=\frac{h}{1-\tau}\, .
  \end{equation}

It turns out that the 
instanton vector bundles defined above with the help of idempotents
$p_n$, $n\in\bZ$, can also be constructed by 
viewing $S^4$ as the quaternionic projective space
$\H P^1$ and considering 
the tautological quaternionic  line bundle.
The bundle can be described by a projection which,
in the foregoing  three coordinate systems, at a generic point of $S^4$, 
can be respectively written as the following matrix in $M_2(\H)$:\\[-6pt]
 \begin{align}
\label{instantonproj}
\half\left[
  \begin{array}{cc}
  1-\tau & \bar h \cr
  h      & 1+\tau
  \end{array}
  \right] 
&\,=\,\frac{1}{\sqrt{2}}\left[
  \begin{array}{c}
  \sqrt{1-\tau}  \cr
  \frac{h}{\sqrt{1-\tau}}      
  \end{array}
  \right]
\left(\frac{1}{\sqrt{2}}\left[
  \begin{array}{c}
  \sqrt{1-\tau}  \cr
  \frac{h}{\sqrt{1-\tau}}      
  \end{array}
  \right]\right)^*,\\[8pt]%
 \left[
  \begin{array}{cc}
  1-t     & \sqrt{t(1-t)}\bar g \cr
 \sqrt{t(1-t)}g & t
  \end{array}
  \right] 
&\,=\,\left[
  \begin{array}{c}
  \sqrt{(1-t)} \cr
  \sqrt{t}g 
  \end{array}
  \right]\left(\left[
  \begin{array}{c}
  \sqrt{(1-t)} \cr
  \sqrt{t}g 
  \end{array}
  \right]\right)^*,\label{gproj}\\[8pt]
  \frac{1}{1+z\bar z}\left[
  \begin{array}{cc}
  1      & \bar z \cr
  z      & z\bar z
  \end{array}
  \right] 
&\,=\,\frac{1}{\sqrt{1+z\bar z}}\left[
  \begin{array}{cc}
  1      \cr
  z      
  \end{array}
  \right]\left(\frac{1}{\sqrt{1+z\bar z}}\left[
  \begin{array}{cc}
  1      \cr
  z      
  \end{array}
  \right]\right)^*,\label{zproj}
  \\[-12pt] \nonumber
  \end{align}
where $\;\bar{}\;$ denotes the quaternionic conjugation and $^*$ stands for the 
quaternionic-Hermitian adjoint. In what follows,
we shall use the middle matrix~\eqref{gproj}. The aforementioned
projection is given by varying this matrix over~$S^4$
viewed as the non-reduced suspension of $S^3$ (see~\eqref{cong}):
\begin{gather}
\label{pn}
S^4
\ni [(t,g)]  \,\stackrel{P_n\;}{\longmapsto}\,  
\left[
  \begin{array}{cc}
  1-t     & \sqrt{t(1-t)}\bar g^n \cr
 \sqrt{t(1-t)}g^n & t
  \end{array}
  \right] 
\in M_{4}(C(S^4)) ,\quad\forall\; n\in \bZ.
 \end{gather}
 \begin{theorem}\label{homotopy}
Let $P_n$, $n\in\bZ$, be the projection \eqref{pn}, 
and let $p_n$, $n\in\bZ$, 
 be an idempotent
corresponding by the formula~\eqref{P}
to $v^n$ for some $v\in M_N(C(S^3))$. 
Then, 
\[\nonumber
K_1(C(S^3))\ni [u]=[v^{-1}] \quad\Longrightarrow\quad
 \forall\; n\in \bZ:\; [P_n]=[p_n]\in K_0(C(S^4)).
\]
 \end{theorem}
 \begin{proof} 
 We prove the theorem by defining a homotopy between $P_n$ and $p_n$. 
The latter is given by \eqref{xysusp} with $a$ replaced by $v^n$.
Since $K_1(C(S^3))\ni [u]=[v^{-1}]$ by assumption, the formula
\eqref{xysusp} with $a$ substituted by $u^{-n}$ yields an idempotent
belonging to the class $[p_n]\in K_0(C(S^4))$.
 Next, we  specify  $\psi \in C([\mbox{$\frac{1}{2}$}, 1])$ 
appearing in \eqref{xysusp} to 
 be a $[0,1]$-valued function 
 such that $\psi(\mbox{$\frac{1}{2}$}) = 1$ and $\psi(1) = 0$. 
 Then we construct the following families of functions
 \[  \label{psis}
 \psi_s(t):= 
\psi\big(\mbox{$\frac{t +1 -s}{2-s}$}\big), 
\quad t\in[\mbox{$\frac{s}{2}$}, 1],\quad s\in [0,1], 
 \]
\begin{align}
X^s_0(t)&:= \left\{
\begin{array}{ccl}
1 &\text{for} & 0\leq t <\mbox{$\frac{s}{2}$},\\
0 &\text{for}& \mbox{$\frac{s}{2}$}\leq t\leq 1,
\end{array}
\right.  &
X^s_1(t)&:= \left\{
\begin{array}{ccl}
0 &\text{for} & 0\leq t <\mbox{$\frac{s}{2}$},\\
\psi_s(t) \big(2-\psi_s(t)^2\big) &\text{for} & \mbox{$\frac{s}{2}$}\leq t\leq 1,
\end{array}
\right. \nonumber\\[4pt]
Y^s_1(t)&:= \left\{
\begin{array}{ccl}
0 &\text{for} & 0\leq t <\mbox{$\frac{s}{2}$},\\
\psi_s(t)  &\text{for}& \mbox{$\frac{s}{2}$}\leq t\leq 1,
\end{array}
\right. 
& X^s_2(t)&:= \left\{
\begin{array}{ccl}
0 &\text{for} & 0\leq t <\mbox{$\frac{s}{2}$},\\
 1-\psi_s(t)^2 &\text{for} & \mbox{$\frac{s}{2}$}\leq t\leq 1,
\end{array}
\right. 
\end{align}
 and define  the family of matrices resembling the matrices
from \eqref{xysusp}
\[
X^s:=\left[
  \begin{array}{c}
  X^s_0\, I_2 +X^s_1\, u^{-n} \cr
  X^s_2\, I_2
  \end{array}
  \right],\quad 
  Y^s:=\left[
X^s_0\, I_2 +Y^s_1\, u^{n} \,,\, X^s_2\, I_2
  \right],
\]
where $I_2$ denotes the $2\times 2$ identity matrix. The
point of this construction is that $X^s \, Y^s$ is an idempotent
matrix in $M_4(\Sigma C(S^3))$.
Indeed, one easily checks that $Y^s \, X^s= I_2$, whence 
$(X^s \, Y^s)^2= X^s \, Y^s$. Furthermore, by a straightforward
computation one verifies that the entries of the matrix $X^s \, Y^s$
are continuous functions from $[0,1]$ to $C(S^3)$ taking numerical values
at the endpoints. Now, 
applying the isomorphism \eqref{via} and comparing with \eqref{xysusp}, 
we see that $X^1\,Y^1\cong p_n$. Hence $p_n$ is isomorphic to an
idempotent homotopic to $X^0\,Y^0$.

On the other hand, let us observe that $P_n=G G^*$, where 
\[
G:=\left[
  \begin{array}{c}
  G_1\, u^{-n} \cr
  G_2\, I_2
  \end{array}
  \right],\quad G_1(t):=\sqrt{1-t},\quad G_2(t):=\sqrt{t},\quad
t\in [0,1]. 
\]
Next, we  define another family of matrices
\[
V^r:=\sqrt{r}\, X^0 +\sqrt{1-r} \,G,\quad W^r:=\sqrt{r}\, Y^0 +\sqrt{1-r}\, G^*,\quad  r\in [0,1]. 
\]
Multiplying these matrices yields $W^r V^r= f^r I_2$, where
\[
f^r(t) := 1+\sqrt{r-r^2}\Big(G_1(t)\big(X^0_1(t)+Y^0_1(t)\big)+ 2 G_2(t) X^0_2(t)\Big),\quad t\in [0,1].
\]
Since $\psi([\half,1])\subseteq [0,1]$, 
one can easily check that  $f^r(t) \geq 1$, for all $t\in [0,1]$. 
Hence we can define
$Q^r := \mbox{$\frac{1}{f^r}$} V^r W^r$. 
As before, by a straightforward computation we verify that
that the entries of each matrix $Q^r$ belong to $\Sigma C(S^3)$. 
Furthermore, all these matrices are evidently idempotent:
\[
(Q^r)^2=
\mbox{$\frac{1}{(f^r)^2}$}V^r W^rV^r W^r=Q^r,\quad r\in [0,1].
\]
Combining it with the fact that  $Q^0= GG^*= P_n$ and $Q^1= X^0\,Y^0$,
we conclude that $P_n$ is homotopic to $X^0\,Y^0$. Since the
latter is homotopic to $X^1\,Y^1\cong p_n$ by the first
part of the proof, by the  homotopy invariance of the K-groups, 
we infer that $[P_n]=[p_n]$ in $K_0(C(S^4))$, as claimed.
\end{proof} 

To end with, let us note that there is yet another way to  describe
the
idempotents $P_n$, $n\in\bZ$. We view $S^4$ as the non-reduced
suspension of $S^3$ via \eqref{cong}, and take advantage of the covering
of $S^4$ by two open balls $U_1:=S^4\setminus\{[(1,g)]\}$ and 
$U_2:=S^4\setminus\{[(0,g)]\}$. Now we can 
take two continuous functions on $S^4$ 
given by $f_1([(t,g)]):=1-t$ and 
$f_2([(t,g)]):=t$, 
and view them as  a kind of
a partition of unity  subordinated to this covering, i.e.\
$f_1+f_2=1$ and $f_i|_{S^4\setminus U_i}=0$, $i\in\{1,2\}$. 
Next, we take the following family of continuous 
functions
\begin{eqnarray}\begin{split}
&\phi_{ij}^n\colon U_i\cap U_j\longrightarrow GL_2(\bC),\quad 
i,j\in\{1,2\}, &
\\
&\phi_{11}^n([(t,g)]):=1=:\phi_{22}^n([(t,g)]),\quad 
\phi_{21}^n([(t,g)]):=u(g)^n,\quad \phi_{12}^n([(t,g)]):=u(g)^{-n},&
\end{split}
\end{eqnarray}
as transition
functions  corresponding to the same covering. Their supports
are not in $U_i\cap U_j$, $i,j\in\{1,2\}$, so that their extensions
 by zero to  functions
on $S^4$ is no longer continuous. However, inserting them into
Karoubi's formula \cite[p.~35]{Kar} together with
$f_1$ and $f_2$ in place of
a partition of unity
subordinate to the covering $\{U_i\}_{i=1}^2$,
we obtain continuous functions on $S^4$ that are 
precisely the entries of the idempotent matrix~$P_n$:
\be\label{steen}
P^{ij}_n:= \sqrt{f_i} \,\phi_{ij}^n\, \sqrt{f_j}, \quad i,j=1,2.
\ee

\subsection{Instanton idempotents from   
Woronowicz's \boldmath$SU_q(2)$}
\label{section_qS4}
 
We take $A:=C(S^3_q)$ to be
the universal unital C*-algebra 
generated by $\a$ and $\b$ subject to the relations
\begin{equation}
\label{suq2_relations}
\a \b = q \b \a ,\,  \;\ \a \b^* = \bar q \b^* \a ,\,  \;\ \b \b^* = \b^* \b,\, \;\
\a^* \a + \b \b^* = 1, \,  \;\ \a \a^* + |q|^2 \b \b^* = 1, \, 
\end{equation}
where $q\in \C$ and $|q|\leq 1$. 
Since $C(S^3_q)\cong C(S^3_{q'})$ for any 
$q$ and $q'$ such that $|q|,|q'|<1$, without the loss of generality
we can restrict to the real $0<q<1$ and unitary $q\in U(1)$ 
cases\footnote{We are grateful to Piotr M.~So\l tan 
for explaining this to us.} (cf.~\cite{garden}).
The case $q =1$ is the classical case discussed earlier,
and the case $q \in U(1)\setminus \{1\}$ is left for the subsequent
section. Herein we focus on the case $0<q<1$, so that $C(S^3_q)$
is the 
C*-algebra of  Woronowicz's quantum group $SU_q(2)$~\cite{wo}.
There is a  faithful representation $\pi$ of $C(S_q^3)$
on the Hilbert space $\ell^2 ( \bN \times \bZ)$ with 
an orthonormal basis 
$\{ e_{m,n} \}_{m\in\bN,\,n\in\bZ}$
given by
\begin{align}
&&\pi(\a) e_{m,n} := \hl_m e_{m-1,n}\,, &&&
\quad \pi(\a^\ast )e_{m,n} := \hl_{m+1} e_{m+1,n}\,,&&\nonumber\\[-8pt]
&& &&& && \label{pialpha} \\[-8pt]
&&\pi(\b) e_{m,n} := q^m e_{m,n+1}\,, &&&
\quad \pi(\b^\ast )e_{m,n} := 
q^m e_{m,n-1} \,,  &&\nonumber
\end{align}
where  $\hl_m := (1 - q^{2m} )^{1/2}$ (see~\cite[Corollary~2.3]{mnw}).

By \cite[Theorem~2.4]{mnw}, 
the K-groups of $C(S_q^3)$ satisfy the assumptions 
of Section~\ref{section_suspension}. Therefore, 
all we need to
construct instanton idempotents  is to find a generator
of $K_1(C(S^3_q))$. Guided by the classical case, we presume that
the $K_1$-class of the fundamental representation of $SU_q(2)$
\begin{equation}
\label{K1gen}
u_q := \left[ 
\begin{array}{cc}
\a & -\bar q \b^\ast \cr
\b & \a^\ast 
\end{array}
\right] \in U_2 (C(S_q^3)) 
\end{equation}
is a desired generator.
It is stated in \cite{Co04}  that the class 
$[u_q]$ is non-trivial in $K_1(C(S_q^3))$,
and its pairing with the $K$-homology class of the spectral 
triple of  \cite{CP} 
is left as an exercise. In \cite{dlsvsv2}, the 
 index pairing of $u_q$ with 
the unbounded $K$-cycle (3-summable spectral triple) 
constructed in \cite{dlsvsv1} 
 was computed  to be~1 (cf.~\cite{lpr}). From the integrality of the index pairing,
 we conclude that   
$[u_q]$ generates $K_1(C(S^3_q))  \cong \bZ$. 
Thus we have shown that
we can use the fundamental representation of $SU_q(2)$ to
construct instanton idempotents of arbitrary charges.
This goes along the lines of the pioneering Pflaum's construction
\cite{pflaum} that was reformulated and continued in 
 \cite{dlm} and \cite{dl}. 

Since finding a $K_1$-generator is pivotal in our construction
of instanton idempotents, we devote the remainder of this section
to two computations showing that a given unitary represents a
$K_1$-generator. As explained above, to prove that the fundamental
representation $u_q$ is such a representative, it suffices to
show that its index pairing with a K-homology class equals to~1.
In a simple and explicit way, we calculate the pairing of $u_q$
with the bounded K-cycle (odd Fredholm module) 
given by \cite[p.~176]{mnw}
\begin{equation}
\label{z}
{\bf z} := ( \ell^2 ( \bN \times \bZ), \pi, F) \,,
\ \; {\rm where}\ \; F e_{m,n} :=
\left\{
\begin{array}{ccl}
\phantom{-}e_{m,n} &\text{for} & n >0,\cr
-e_{m,n}  &\text{for} & n \leq 0.
\end{array}
\right.
\end{equation}
Note first that it follows from the definition \eqref{pialpha} of $\pi$ 
 that  
$[ {F} , \pi(\a)]  = 0= [ {F} , \pi(\a^\ast) ]$ and 
\begin{eqnarray}
&
\pi(\b^\ast) [F,\pi(\b)] e_{m,n} &\,\,= \,
\left\{
\begin{array}{ccl}
\phantom{\mbox{$-$}}2 \,q^{2m } \,  e_{m, 0} &\text{for} &n = 0,\cr
0 &\text{for} & n \neq 0,
\end{array}
\right.\\
&\pi(\b^\ast) [F,\pi(\b)] e_{m,n} &\,\,= \,
\left\{
\begin{array}{ccl}
 - 2 \,q^{2m} \, e_{m, 1} &\text{for} & n = 1,\cr
0 &\text{for} & n \neq 1.
\end{array}
\right.
\end{eqnarray}
Furthermore, observe that
 all commutators of $F$ with the generators of $\pi(C(S_q^3))$ 
are of trace class. 
Therefore, ${\bf z}$ is 1-summable 
over the *-algebra generated by $\a$ and $\b$. 
Next, let $\Tr$ denote the operator trace on $\ell^2( \bN \times \bZ)$ 
and $\Tr_2$ the usual trace on $2\times 2$-matrices. Since $F^2=\id$,
for
a trace-class operator $T$ we immediately infer that 
$\Tr(F(FT+TF))=2\Tr(T)$. Combining this equality with the 1-summability
of $\bf z$,
by
\cite[Definition IV.1.3, Propositions~III.3.3 and IV.1.4]{connes_book},
the index pairing takes the following form: 
\begin{align}\nonumber
<[{\bf z}],[u_q]> &=
\mbox{$\half$}\Tr\left(\Tr_2\left( 
\left[ 
\begin{array}{cc}
\pi(\a^\ast) -1 &  \pi(\b^\ast )\cr
- q \pi(\b) & \pi(\a)-1
\end{array}
\right]
\left[ 
\begin{array}{cr}
[F, \pi(\a)] & -q [F, \pi(\b^\ast) ]\cr
[F, \pi(\b)] & [F, \pi(\a^\ast)] 
\end{array}
\right] \right)\right)\nonumber\\
&=\mbox{$\half$}\Tr \left( \pi(\b^\ast) [F,\pi(\b)] -q^2 \pi(\b^\ast) [F,\pi(\b)]
\right) \nonumber  \\
&= (1-q^2) \sum_{m=0}^\infty q^{2m}\nonumber\\
&=1.
\end{align}

On the other hand, without taking advantage of the index pairing,
it was already proven  in \cite[Theorem~2.4]{mnw}
that a generator of~$K_1(C(S^3_q))$ can be represented by
\[\label{japw}
w:=\b e(0) + (1- e(0)) \in C(S_q^3),
\]
 where
 $e(0)$ is the spectral projection  corresponding
to the isolated point $1$ of the spectrum of $\b^\ast \b$.
The argument used in \cite{mnw} relied on the six-term exact 
sequence in K-theory applied to the short exact sequence of C*-algebras
\[
0\longrightarrow \mathcal{K}\otimes C(S^1)
\longrightarrow C(S^3_q)\longrightarrow C(S^1)
\longrightarrow 0,
\]
where $\mathcal{K}$ denotes the ideal of compact operators
and $C(S^1)$ is the C*-algebra of all continuous functions on
the unit circle.
Herein we apply our pullback point of view also to~$C(S^3_q)$.
We present it as a pullback \mbox{C*-algebra} and compute
an even-to-odd connecting homomorphism of the Mayer-Vietoris
six term exact sequence using the formula~\eqref{baumhajac}. 
Thus we obtain an alternative proof
that $w$ represents a generating class of~$K_1(C(S^3_q))$.

To begin with, recall that it was shown in \cite{HW} that 
$C(S^3_q)$ can be given by the 
following pullback diagram: 
\begin{equation}\label{barP}
\xymatrix{
& \makebox[100pt][c]{$\mbox{ }\hspace{-24pt} 
C(S^3_q)\;\cong \; \CS \underset{{\scriptsize \CS\ot\CS}}{\times}\sT
\ot\,\CS   $}
\ar[dl]_{\pr_1} \ar[dr]^{\pr_2}& \\
\CS\ar[rd]_{\pi_1:=\Delta\ \ } &    &  
\sT \ot\CS \ar[ld]^{\ \pi_2:=\sigma\otimes\id}\\
 &\mbox{ }\hspace{-30pt} C(T^2)=
\CS\ot\CS\,.\phantom{\hat{\hat{I}}}&   }
\end{equation}
Here $\Delta$ is the pullback of the multiplication map restricted
to unitary complex numbers and 
$\sT$ denotes the Toeplitz algebra viewed as the  C*-algebra generated 
by the unilateral shift 
$S$ acting on the Hilbert space $\lN$. The symbol map $\sigma: \sT\rightarrow \CS$ is given 
by $\sigma(S) = U$, where 
$U$ denotes the unitary generator of $\CS$,
and $\mathrm{pr}_i$, $i\in\{1,2\}$, are the restrictions of
canonical projections. 

For an explicit description of the isomorphism given in \cite{HW}, 
we first realize $\CS$ as an operator algebra on 
the Hilbert space $\lZ$ via
the faithful representation $\rho$
 identifying the unitary generator $U$ 
with the negative bilateral shift given on an
orthonormal basis of $\lZ$ by
\[
\rho(U):=V,\quad Ve_k:=e_{k-1}\,,\quad k\in\bZ. 
\]
On the other hand, since $\pi_1:=\Delta$ is injective, 
${\mathrm{pr}_2}$ corestricted to its image is an 
isomorphism. Combining the above identifications, 
we obtain
\begin{equation}  \label{SUq-isom}
C(S^3_q)\, \cong \, P:=\overline{\mathrm{span}}
\{ T\otimes V^N \in  \sT \ot\rho(\CS)\; |
\; \sigma(T)=0 \ \,\text{or}\,\  
\sigma(T)=U^N,\ N\in\bZ\}, 
\end{equation}
where $\overline{\mathrm{span}}$ denotes the closed linear span. 
Thus we view $C(S^3_q)$ as an algebra of bounded operators on the
Hilbert space $\ell^2(\bN\ot\bZ)$.
To define an isomorphism
$\psi\colon C(S^3_q)\to P$ (see~\cite{HW}) on 
the generators from~\eqref{suq2_relations}, 
consider the bounded linear operators  $y$ and $z$ on the
Hilbert space $\lN$  given by 
\[
ye_n :=q^n e_n,\quad z e_n := \lambda_{n+1} e_{n+1}\,,\quad n\in\bN.
\] 
Here $\{e_n\}_{n\in\bN}$ is an orthonormal basis of 
$\lN$ and the $\lambda_n$'s are as in~\eqref{pialpha}.
Furthermore,
recall that the symbol map yields the following
well-known short exact sequence of \mbox{C*-algebras}
\begin{equation}\label{sesToeplitz}
\xymatrix{
 0\,\ar[r]&\,  \cK(\lN)\,\ar@{^{(}->}[r] &\, \sT 
\,\ar[r]^{\sigma\quad} &\, \CS\,\ar[r] &\,0\,.}
\end{equation}
Clearly, $y\in  \cK(\lN)\subset \sT$. 
Also, since $z-S\in \cK(\lN)$, we infer that
 $z\in \sT$ and $\sigma(z)=\sigma(S)=U$. 
Putting it all together, we conclude that
the isomorphism $\psi$ is given by
\begin{align}
&&&\psi(\a) := z^*\otimes V^*, && \psi(\a^\ast )
:= z\otimes V,\nonumber&&\\[-8pt]
&&& && &&  \\[-8pt]
&&&\psi(\b)  := y\otimes V^*, && \psi(\b^\ast ) 
:= y \otimes V.\nonumber&&
\end{align}

With the isomorphism $\psi$ at hand, we can  
describe the unitary $w$ from 
\eqref{japw} in terms of 
 generators of $C(S^3_q)$. Indeed,
using $z|z|^{-1}=S$ and the fact the spectral projection of $y^2$
corresponding to the eigenvalue $1$ equals $1-SS^*$, we get 
\begin{align} \nonumber
\psi( 1+(\hb-1)(1-\ha^*(\ha \ha^*)^{-1}  \ha) ) &= 1\ot 1 + ( y \otimes V^* -1\ot 1) (1\ot 1 - SS^*\ot 1)\\
&= SS^*\ot 1 +     y(1-SS^*)\ot V^*\nonumber\\
&=\psi(1-e(0)) + \psi(\b e(0))\nonumber\\
&=\psi(w) \label{u},
\end{align}
so that
\[  \label{upsi}
1+(\hb-1)(1-\ha^*(\ha \ha^*)^{-1}  \ha)=w \in C(S^3_q).
\] 

We want to prove that $[w]$ is a generator of $K_1(C(S^3_q))$
by showing that it coincides with the even-to-odd connecting homomorphism \eqref{baumhajac}
of the Mayer-Vietoris six-term exact sequence 
of the fibre product \eqref{barP}
applied to a generator of the non-trivial part 
of $K_0(C(T^2))\cong\bZ\oplus\bZ$. (The trivial part is generated
by $1\in C(T^2)$.)
To this end, we follow \cite{L}
to explicitly represent such a generator
by the projection 
\begin{equation}
\label{K0gen}
e := \left[ 
\begin{array}{cc}
f & g+U h \\
g+U^* h& 1-f
\end{array}
\right] \in M_2 (C(T^2)) \,. 
\end{equation}
Here we view functions on $T^2$ as periodic functions on $\bR^2$. Thus the unitary generator of $\CS$
is now a function given by $U(t):=\ex^{2\pi i t}$,  and 
\begin{align} 	\label{f}
f(s)&:= \left\{
\begin{array}{ccl}
1-2s &\text{for} & s\in[0,\mbox{$\frac{1}{2}$}]\ \mathrm{mod}\ 1,\\
1+2s &\text{for} & s\in[-\mbox{$\frac{1}{2}$},0]\ \mathrm{mod}\ 1,
\end{array}
\right.\\[4pt]
h(s)&:= \chi_{[0,\frac{1}{2}]}(s) \sqrt{f(s)(1-f(s))},\\[4pt]
g(s)&:= \big(1-\chi_{[0,\frac{1}{2}]}(s)\big) \sqrt{f(s)(1-f(s))}
\end{align}
with 
\[ \label{chi0}
\chi_{[0,\frac{1}{2}]}(s):= \left\{
\begin{array}{ccl}
1 &\text{for} & s\in[0,\mbox{$\frac{1}{2}$}]\ \mathrm{mod}\ 1,\\
0 &\text{for} & s\in (-\mbox{$\frac{1}{2}$},0)\ \mathrm{mod}\ 1.
\end{array}
\right.
\]
Note that $f,g,h\in\CS$ and $e^*=e$. Furthermore, from 
\[  \label{fgh}
 gh=0,\qquad  f^2+g^2+h^2=f,
\]
it follows that $e^2=e$. As in \cite{L},  by checking that 
$\mathrm{ch_1}(e)=1$, one can prove that 
$[e]$ generates the non-trivial part of~$K_0(C(T^2))$. 

\begin{theorem} 
Let $\ha$ and $\hb$ be the generators of $C(S^3_q)$ given 
in~\eqref{suq2_relations}, and let $e$ be
the projection generating the non-trivial part of the $K_0$-group of the 
2-torus provided in~\eqref{K0gen}. 
Then the even-to-odd connecting homomorphism~\eqref{baumhajac} applied 
to the $K_0$-class of~$e$ yields 
\[
\partial_{01}([e])=[1+(\hb-1)(1-\ha^*(\ha \ha^*)^{-1}  \ha)]
=1\in K_1(C(S^3_q)=\mathbb{Z}\;.
\]
\end{theorem}
\begin{proof}
We begin by showing that $\partial_{01}([e])$  generates 
 $K_1(C(S^3_q)$. 
To this end, we consider the 
Mayer-Vietoris six-term exact sequence of the fibre 
product~\eqref{barP}:
 \begin{equation}  
     	\begin{CD}
	{K_0 (C(S^3_q))} @ >{({\pr_1}_\ast , {\pr_2}_\ast)}>> 
	{K_0 ( \CS ) \oplus K_0 (\sT \ot\CS)} @ >{{\pi_2}_\ast - {\pi_1}_\ast} 
	>> {K_0 ({C(T^2)} )} @ . @  .\\ 
 	@ A{\partial_{10}}AA @ . @ VV
     {\partial_{01}} V @ .\\
	{K_1 ({C(T^2)} ) } @ <{{\pi_2}_\ast - {\pi_1}_\ast} << 
	{K_1 (\CS ) \oplus K_1 (\sT \ot\CS)} @ <{({\pr_1}_\ast , {\pr_2}_\ast)} 
	<< {K_1 (C(S^3_q))\,.} @ . @ .\\
     	\end{CD}
\end{equation}
The K-groups of $C(S^3_q))$ were computed in \cite[Theorem~2.4]{mnw} 
to be $\bZ$, and other K-groups are well-known to be $\bZ$ or 
$\bZ\oplus\bZ$ \cite[p.\ 123]{wegge-olsen}. Since both
$K_0 ( \CS )$ and \mbox{$K_0 (\sT \ot\CS)$} are generated by the 1's of the
C*-algebras, the induced map to $K_0 ({C(T^2)} )$ is given by
$(m,n)\mapsto (m-n,0)$. On the other hand, since 
$K_1 ( \CS )$
is generated by~$[U]$, \mbox{$K_1 (\sT \ot\CS)$} is generated by
$[1\ot U]$, and $K_1 ({C(T^2)} )=K_1 ({C(S^1)\ot C(S^1)} )$
is generated by $[U\ot 1]$ and $[1\ot U]$,
the induced map to $K_1 ({C(T^2)} )$ is given by the formula
$(m,n)\mapsto (-m,n-m)$. Consequently, the exactness of the sequence
forces the preceding map to be zero. Summarizing, we have
arrived at the following exact sequence:
\begin{equation}    \label{ZZ}
     	\begin{CD}
	{\bZ} @ >{\phantom{\,(-m,n-m)\,\mapsfrom\, (m,n)}}>> 
	{\bZ \oplus \bZ} @ >{(n,m)\,\mapsto\, (n-m,0)\ } 
	>> {\bZ\oplus\bZ} @ . @  .\\ 
 	@ A{}AA @ . @ VV
     {\partial_{01}} V @ .\\
	{\bZ\oplus\bZ} @ <{\,(-m,n-m)\,\mapsfrom \,(m,n)} << 
	{\bZ \oplus \bZ} @ <{\phantom{xxxxxx}0\phantom{xxxxxxx}} 
	<< {\,\bZ\,.} @ . @ .\\
     	\end{CD}
\end{equation}
The exactness of the above sequence immediately
implies that $\partial_{01}([e])$ generates the group $K_1 (C(S^3_q))$,
 as claimed.

It remains to compute $\partial_{01}([e])$ in terms of the generators
$\alpha$ and~$\beta$. We do it with the help of~\eqref{baumhajac}.
First, we lift $U\in\CS$ to $S\in\sT$, and then $e\in M_2 (C(T^2))$ to 
\begin{equation}
\label{lift}
Q := \left[ 
\begin{array}{cc}
1\ot f & 1\ot g+S\ot h \\
1\ot g+ S^*\ot h & 1\ot(1-f)
\end{array}
\right] \in M_2 (\sT\ot \CS) \,.
\end{equation}
Next, remembering $S^*S=1$ and \eqref{fgh}, we obtain 
\begin{align} 
\nonumber
Q^2 &= \left[ 
\begin{array}{cc}
1\ot (f^2+g^2) + SS^*\ot h^2 & 1\ot g+S\ot h \\
1\ot g+ S^*\ot h & 1\ot(1-f)
\end{array}
\right] \\
&= Q- 
 \left[ 
\begin{array}{cc}
(1-SS^*)\ot h^2 & 0 \\
0 & 0
\end{array}
\right] .\label{lift3}
\end{align}
Furthermore, observe that for any $\varphi\in\CS$ with 
$\mathrm{supp}(\varphi)\subset [0,\frac{1}{2}]\; \mathrm{mod}\; 1$, 
we have 
\[\label{lift2}
 Q \left[ 
\begin{array}{cc}
(1-SS^*)\ot \varphi & 0 \\
0 & 0
\end{array}
\right] 
=
\left[ 
\begin{array}{cc}
(1-SS^*)\ot  f\varphi& 0 \\
0 & 0
\end{array}
\right]\\
= 
 \left[ 
\begin{array}{cc}
(1-SS^*)\ot \varphi & 0 \\
0 & 0
\end{array}
\right] Q
\]
because $(1-SS^*)S=S^*(1-SS^*)=0$ and $g\varphi=0$. 
For the sake of the induction argument, assume that 
\begin{equation} 
\label{lift4}
Q^n 
= Q- 
 \left[ 
\begin{array}{cc}
(1-SS^*)\ot \chi_{[0,\frac{1}{2}]}(f-f^n) & 0 \\
0 & 0
\end{array}
\right].
\end{equation}
Then combining \eqref{lift3} with $h^2=\chi_{[0,\frac{1}{2}]}(f-f^2)$
and \eqref{lift2}, we compute:
\begin{align}
Q^{n+1} 
&= Q^2- 
 \left[ 
\begin{array}{cc}
(1-SS^*)\ot \chi_{[0,\frac{1}{2}]}(f-f^n) & 0 \\
0 & 0
\end{array}
\right]Q \nonumber\\
&= Q
-\left[ 
\begin{array}{cc}
(1-SS^*)\ot \chi_{[0,\frac{1}{2}]}(f-f^2) & 0  \\
0 & 0
\end{array}
\right]
- \left[ 
\begin{array}{cc}
(1-SS^*)\ot \chi_{[0,\frac{1}{2}]}f(f-f^n) & 0 \\
0 & 0
\end{array}
\right]\nonumber \\
&= Q- 
 \left[ 
\begin{array}{cc}
(1-SS^*)\ot \chi_{[0,\frac{1}{2}]}(f-f^{n+1}) & 0 \\
0 & 0
\end{array}
\right].
\end{align}
By induction, this proves the formula~\eqref{lift4} for $Q^n$
for all $n\in\bN$.
Hence, noting the equality
$\chi_{[0,\frac{1}{2}]}(\exp(2\pi i f)-1)= \exp(2\pi i 
\chi_{[0,\frac{1}{2}]}f)-1$,  we can compute: 
\begin{align}
\exp(2\pi i Q)
&=\left[ 
\begin{array}{cc}
1\ot 1 & 0 \\
0 & 1\ot 1
\end{array}
\right] 
+ \sum_{n=1}^{\infty}\frac{(2\pi i)^n}{n!} 
\left(Q- 
 \left[ 
\begin{array}{cc}
(1-SS^*)\ot \chi_{[0,\frac{1}{2}]}(f-f^n) & 0 \\
0 & 0
\end{array}
\right] 
\right)\nonumber\\
&=  \left[ 
\begin{array}{cc}
1\ot 1 & 0 \\
0 & 1\ot 1
\end{array}
\right] +\sum_{n=0}^{\infty}\frac{(2\pi i)^n}{n!} 
\left(Q- 
 \left[ 
\begin{array}{cc}
(1-SS^*)\ot \chi_{[0,\frac{1}{2}]}(f-f^n) & 0 \\
0 & 0
\end{array}
\right] 
\right) \nonumber\\
&\phantom{=}\ - \left(Q- 
 \left[ 
\begin{array}{cc}
(1-SS^*)\ot \chi_{[0,\frac{1}{2}]}(f-1) & 0 \\
0 & 0
\end{array}
\right] 
\right)\nonumber\\
&= \left[ 
\begin{array}{cc}
1\ot 1 + \big(1-SS^*\big)\ot  
\big(\exp(2\pi i \chi_{[0,\frac{1}{2}]}f)-1\big) 
& 0 \\
0 & 1\ot 1
\end{array}
\right] .
\end{align}
Now, by inserting the definition \eqref{f} of $f$, 
remembering that the pullback C*-algebra
and $P$ defined in~\eqref{SUq-isom} are isomorphic via
the corestriction of 
\mbox{$(\id\ot\rho)\circ\pr_2$} 
to its image, and removing the lower diagonal 
entry $1\ot 1$ 
by the defining equivalence relation of the $K_1$-group,
we can write 
\begin{align}\nonumber
\partial_{01}([e])&=[((\id\ot\rho)\circ\pr_2)(I,\exp(2\pi i Q)) ] \\
&= [1\ot 1 + \big(1-SS^*\big)\ot  
\rho\big(\exp(-4\pi i\chi_{[0,\frac{1}{2}]} \id_\bR)-1\big)]\,. 
\label{lift5}
\end{align}
Since $1-SS^*$ is a projection, one easily checks that 
$v:=1\ot 1 + \big(1-SS^*\big)\ot  \big(\nu -1\big)$ is unitary for 
any unitary element $\nu\in \CS$. Furthermore, 
$v$ 
belongs to $C(S^3_q)$ by \eqref{SUq-isom} because $\sigma (1-SS^*)=0$. 
Therefore we can define a homotopy of unitaries in $C(S^3_q)$ by 
\[ \label{ut}
v_t:= 1\ot 1 + \big(1-SS^*\big)\ot  
\rho\left(\exp\left(-\mbox{$\frac{4\pi i \chi_t\id_\bR}{1+t}$}\right)-1
\right),\quad t\in[0,1],
\]
where
\[
\chi_t(s):= \left\{
\begin{array}{ccl}
1 &\text{for} & s\in[0,\mbox{$\frac{1}{2}(1+t)$}]\ \mathrm{mod}\ 1,\\
0 &\text{for} & s\notin[0,\mbox{$\frac{1}{2}(1+t)$}]\ \mathrm{mod}\ 1.
\end{array}
\right.
\]
Hence, by the homotopy invariance of K-theory, it follows 
from~\eqref{lift5} that 
$\partial_{01}([e])=[v_0]=[v_1]$ in $K_1(C(S^3_q)$.
Now it suffices to note that $\rho(\exp(-2\pi i\,\id_\bR))=\rho(U^*)
=V^*$  to conclude that 
\[
\partial_{01}([e]) = [v_1]= [1\ot 1 + (1-SS^*)\ot  (V^*-1)]= [SS^*\ot 1 
+ (1-SS^*)\ot  V^*].
\] 
Combining it with the fact that $1-SS^*=y(1-SS^*)$ and \eqref{u} ends the proof.
\end{proof}

\subsection{Instanton projections from Matsumoto's 
 \boldmath$S^3_\theta$}
\label{section_thetaS4}

In this section, we take our deformation parameter $q$ to be unitary, and write 
it $q:= \ex^{2\pi i\theta}$,
$\theta\in [0,1]$. Thus our main building block is the C*-algebra 
$A:= C( S^3_\theta) $ of the Matsumoto quantum 3-sphere~\cite{matsumoto}.
It is defined as the universal unital C*-algebra generated by $\a$, $\b$,
satisfying the relations \eqref{suq2_relations} for $q\in U(1)$.
By \cite[Proposition~2.1.(i)]{mt92} it is a nuclear C*-algebra, so that
\eqref{compatibility} holds for its non-reduced suspension pullback 
diagram.

As in the previous case,
 the K-groups of $C(S^3_\theta)$ satisfy  the assumptions of 
Section~\ref{section_suspension}. We know that the unitary
 matrix
\begin{equation}
\label{K1genn}
\mathrm{u}_\theta := \left[ 
\begin{array}{cc}
\a & -\bar\lambda \b^\ast \cr
\b & \a^\ast 
\end{array}
\right] \in U_2 (C( S^3_\theta) )
\end{equation}
generates $K_1(C( S^3_\theta))$
because
its pairing with a certain spectral triple (unbounded Fredholm module) 
over $C(S^3_\theta)$ is~1. This was  shown
 in \cite{notes} using a Dirac operator isospectral to the
standard classical Dirac operator on $S^3$. Therefore,  by plugging in
$\mathrm{u}_\theta^n$, $n\in\bZ$,  for $a$ in~\eqref{xysusp} and applying \eqref{main2}, 
we conclude that $XY$ is a
$\theta$-instanton idempotent 
 of  charge~$n$ with entries in
$C(S^4_\theta) := \Sigma C(S^3_\theta )$. Finally, by a general argument,
we can always choose a projection in the $K_0$-class of any of our
 idempotents.

This gives a significant simplification of a more direct argument
combining \cite[Proposition~4.14]{BL}
with the computation of the  index pairing of the basic ($n=1$)
$\theta$-instanton projection with an appropriate spectral triple
over $C(S^4_\theta)$. The aforementioned index pairing
 was calculated in \cite{cl01} by applying the local index theorem of 
Connes and Moscovici,
which required
direct computation
containing a few hundreds terms. This determined the charge of the
basic $\theta$-instanton to be equal to 1. With this result at hand,
$\theta$-instanton projections of arbitrary charge were produced 
in~\cite{BL}.

\subsection{Other instanton idempotents}
\label{comments}

The two cases of quantum deformations of $S^3$ 
presented in the previous sections yield
 different types of  C*-algebras. However, they find a common
 denominator in the C*-algebra of the Heegaard quantum $S^3$
\cite{bhms06}.
The even-to-odd connecting homomorphism \eqref{baumhajac} applied to the
Powers-Rieffel idempotent \cite{R} yields a unitary whose powers
allow one to construct a new type of noncommutative instantons~\cite{dhw}.

\subsection*{Acknowledgments}

It is a pleasure to thank N.~Higson,
U.~Kr\"ahmer, G.~Landi,\linebreak R.~Meyer, R.~Nest, P.~M.~So\l tan
 and S.~L.~Woronowicz 
for very useful discussions. 
T.~H.\ and R.~M.\ thank Katedra Metod Matematycznych Fizyki, 
Uniwersytet Warszawski, 
for its support via the EU Transfer of Knowledge contract 
MKTD-CT-2004-509794. 
T.~H.\ 
 is also very grateful to Instytut Matematyczny, Polska Akademia Nauk, 
and to the Isaac Newton Institute, Cambridge, for their support and 
hospitality.


\begin{thebibliography}{88}

\bibitem{atiyah} M.~F.~Atiyah, 
\emph{The Geometry of Yang-Mills Fields}, 
Fermi Lectures, Scuola Normale Superiore Pisa, Pisa, 1979. 

\bibitem{bhms06}  
P.~F.~Baum, P.~M.~Hajac, R.~Matthes and W.~Szyma\'nski, 
\emph{The K-theory of Heegaard-type quantum 3-spheres},
$K$-Theory 35 (2005), 159--186.

\bibitem{bhr}  P.~F.~Baum, P.~M.~Hajac and J.~Rudnik,
\emph{The K-theory of quantum real projective planes},
in preparation.

\bibitem{bm} 
P.~F.~Baum and R.~Meyer,
\emph{The Baum--Connes conjecture, localisation of categories, and quantum groups}, 
in: P.~M.\ Hajac (ed.), 
Lecture Notes on Noncommutative Geometry and Quantum Groups, 
EMS Publ. House, to appear. 

\bibitem{blackadar}
B.~Blackadar, \emph{K-Theory for Operator Algebras}, 
MSRI publications  5, Cambridge University Press, New York, 1998.

 \bibitem{BCT1} 
F.~Bonechi, N.~Ciccoli and M.~Tarlini,
\emph{Noncommutative instantons on the 4-sphere from quantum groups},
Commun.\ Math.\ Phys.\  226 (2002), 419--432.
 
\bibitem{BL} 
S.~Brain and G.~Landi,
\emph{Families of monads and instantons from a noncommutative ADHM 
construction},
 Clay Math.\ Proc.\ 11 (2010),  55--84.  

\bibitem{CP}    
P.~S.~Chakraborty and A.~Pal, 
\emph{Equivariant spectral triples on the quantum SU(2) group}, 
$K$-Theory 28 (2003), 107--126.

\bibitem{connes_book} A.~Connes, 
\emph{Noncommutative Geometry}, 
Academic Press, San Diego, 1994. 

\bibitem{Co04}
A.~Connes,
\emph{Cyclic cohomology, quantum group symmetries and the local index 
formula for $SU_q(2)$},
J.\ Inst.\ Math.\ Jussieu  3 (2004), 17--68.

\bibitem{cl01} A.~Connes and G.~Landi, 
\emph{Noncommutative manifolds, the instanton algebra and isospectral 
deformations}, 
Commun.\ Math.\ Phys.\ 221 (2001), 141--159.

\bibitem{garden} L.~D\k abrowski,
\emph{The garden of quantum spheres},
Banach Center Publications  61 (2003), 37--48.

 \bibitem{dhw}    
L.~D\k abrowski, P.~M.~Hajac, and E.~Wagner,
\emph{From Heegaard unitaries to instanton projections},
in preparation. 

\bibitem{dl}
L.~D\k abrowski and G.~Landi, 
\emph{Instanton algebras and quantum 4-spheres},  
Differential Geom.\ Appl.\  16 (2002), 277--284.

 \bibitem{dlm}
L.~D\k abrowski, G.~Landi and T.~Masuda,
\emph{Instantons on the quantum 4-Spheres $S^4_q$},
Commun.\ Math.\ Phys.\  221 (2001), 161--168.

\bibitem{dlsvsv1}  
L.~D\k abrowski, G.~Landi, A.~Sitarz,  W.~van Suijlekom and J.~C.~Varilly,  
\emph{The Dirac operator on $SU_q(2)$},  
Commun.\ Math.\ Phys.\  259 (2005), 729--759.

\bibitem{dlsvsv2}  
L.~D\k abrowski, G.~Landi, A.~Sitarz, W.~van Suijlekom and J.~C.~Varilly, 
\emph{The local index formula for $SU_q(2)$},
$K$-Theory 35 (2005), 375--394. 

\bibitem{hrz}
P.~M.~Hajac, A.~Rennie and B.~Zieli\'nski,
\emph{The K-theory of Heegaard quantum lens spaces}, 
to appear in J.~Noncommutat.~Geom. 

\bibitem{HW}
P.~M.~Hajac and E.~Wagner,
\emph{The pullbacks of principal coactions}, 
arXiv:1001.0075. 

\bibitem{Kar}
M.~Karoubi,
\emph{K-Theory: An Introduction},
Springer, Berlin, 1978.

\bibitem{notes}
G.~Landi, unpublished notes. 

\bibitem{lpr} 
G.~Landi, C.~Pagani and C.~Reina, 
\emph{A Hopf bundle over a quantum four-sphere from the symplectic group},
Commun.\ Math.\ Phys.\  263 (2006), 65--88.

\bibitem{L} 
T.~A.~Loring, 
\emph{The torus and noncommutative topology}, Ph.D.-thesis, 
University of California, Berkeley, 1986. 

\bibitem{mnw} 
T.~Masuda, Y.~Nakagami and J.~Watanabe, 
\emph{Noncommutative differential geometry on the quantum $SU(2)$, I: An algebraic viewpoint}, 
$K$-Theory 4 (1990), 157--180.

\bibitem{matsumoto} 
K.~Matsumoto, 
\emph{Noncommutative three-dimensional spheres},
Japan.\ J.\ Math\. (N.~S.) 17 (1991), 333--356.

\bibitem{mt92} 
K.~Matsumoto, J.~Tomiyama,
\emph{Non-commutative lens spaces},
J.\ Math.\ Soc.\ Japan 44 (1992), 13--41.

\bibitem{milnor} 
J.~Milnor,
\emph{Introduction to Algebraic $K$-Theory}, Annals of Mathematics Studies 72, 
Princeton University Press, Princeton, New Jersey, 1971.

\bibitem{Pedersen} G.~K.~Pedersen, 
\emph{C*-Algebras and their Automorphism Groups},      
Academic Press, London, 1979.

\bibitem{pflaum} 
M.~J.~Pflaum,
\emph{Quantum groups on fibre bundles},
Commun.\ Math.\ Phys.\ 166 (1994), 279--316.

\bibitem{Ran} 
A.~Ranicki, 
\emph{An explicit projection}, 
http://www.maths.ed.ac.uk/$\sim$aar/surgery/explicit.pdf. 

\bibitem{R} 
M.~A.~Rieffel,
\emph{C*-algebras associated with irrational rotations}, 
Pacific J.\ Math.\ 93 (1981), 415--429.

\bibitem{r-ma90} 
M.~A.~Rieffel,
\emph{Non-commutative tori --- A case study of non-commutative 
differentiable manifolds}, Contemporary Mathematics 105 (1990), 191--211.

\bibitem{rs87} 
J.~Rosenberg and  C.~Schochet, 
\emph{The K\"unneth theorem and the universal coefficient theorem for 
Kasparov's generalized $K$-functor},  
Duke Math.\ J.\  55  (1987), 431--474.

\bibitem{wegge-olsen}
N.~E.~Wegge-Olsen, 
\emph{K-Theory and C*-Algebras: A Friendly Approach}, 
Oxford University Press Inc., New York, 1993.

\bibitem{wo} 
S.~L.~Woronowicz, 
\emph{Twisted $SU(2)$ group: an example of a noncommutative differential calculus},
Publ.\ RIMS, Kyoto Univ. 23 (1987), 117--181.

\end{thebibliography}
\end{document}